\documentclass[11pt]{amsart}
\usepackage{amssymb,amsfonts,amscd,verbatim}
\usepackage[leqno]{amsmath}
\usepackage{mathrsfs}
\usepackage{xy} % for diagrams

 \swapnumbers

\numberwithin{equation}{subsection}%% eq. are numbered as ?.?.?
\newtheorem{thm}[equation]{Theorem} % thm and equation share the
\newtheorem{propose}[equation]{Proposition}
\newtheorem{lemma}[equation]{Lemma}

\theoremstyle{definition}
\newtheorem{defn}[equation]{Definition}

\newtheorem{examples}[equation]{Examples}

\newcommand{\coker}{\operatorname{coker}} % For the cokernel of a morphism
\newcommand{\Shv}{{\rm Shv}}

\newcommand{\Ext}{{\rm Ext}}
\newcommand{\Hom}{{\rm Hom}}
\newcommand{\Spec}{{\rm Spec}}
\newcommand{\by}[1]{\stackrel{#1}{\rightarrow}}
\newcommand{\nobd}{\nobreakdash}

\newcommand{\G}{\mathbb{G}}     % For group schemes
\newcommand{\C}{\mathbb{C}}     % For complex numbers
\newcommand{\Z}{\mathbb{Z}}     % For integers

\newcommand{\df}{\mbox{\,${:=}$}\,}

\newcommand{\ie}{{i.e.},\ }
\newcommand{\cf}{{\it cf.\/}\ }

\newcommand{\id}{{\rm id}}
\newcommand{\eff}{{\rm eff}}
\newcommand{\et}{{ \rm {\acute et}}}
\newcommand{\tor}{{\rm tor}}

\newcommand{\fr}{{\rm fr}}

\newcommand{\fppf}{{\rm fppf}}
\newcommand{\M}{\mathcal{M}_1}   % For category of Deligne 1-Motives
\newcommand{\tM}{{}^t\!\mathcal{M}_1} % For category of 1-Motives
\newcommand{\Ma}{ \mathcal{M}_1^{\rm a}}   % For category of Laumon's 1-Motives
\newcommand{\tMa}{ {}^t\!\mathcal{M}_1^{\rm a}} % For category of 1-Motives
\newcommand{\tMaeff}{{}^t\! \mathcal{M}_1^{\rm a,eff} }   % For category of effective Laumon's 1-Motives

\begin{document}

\xyoption{all}%%

\title{Generalized $1$-motivic sheaves}
\author{A. Bertapelle}
\address{Dipartimento di Matematica, Universit\`a degli Studi di Padova\\Via Trieste, 63\\Padova-- I-35121\\ Italy}
\email{alessandra.bertapelle@unipd.it}
\thanks{Partially supported by   Progetto di Eccellenza Cariparo 2008-2009
``Differential methods in Arithmetics, Geometry and Algebra''  }
\subjclass[2010]{14C15,18E30,14F42,14L15}

\begin{abstract} {\normalsize We extend the construction of the category of $1$-motivic sheaves in \cite{BVK} allowing  quotients of connected algebraic $k$-groups by formal $k$-groups. We show that  its bounded derived category is equivalent to the bounded derived category of the category of generalized $1$-motives with torsion introduced in \cite{BVB}.} \end{abstract}

\maketitle

%%%
\section{Introduction}\label{sec.intro}

Let $k$ be a field of characteristic $0$. There are several generalizations of the category of Deligne $1$-motives, \ie the category of complexes $[L\to G]$, in degrees $-1,0$,  where $G$ is  a semiabelian $k$-variety and $L$ is an \'etale group scheme over $k$, isomorphic to a $\Z^r$ over the algebraic closure of $k$. 
One generalization is due to Laumon. In \cite{LAU} he considers as $G$ any connected (commutative) algebraic $k$-group and allows copies of the formal $k$-group $\hat \G_a$  in degree $-1$. Another generalization, due to Barbieri-Viale, Rosenschon and Saito (see \cite{BVRS}), allows torsion \'etale groups in degree $-1$ and, up to localization at the class of quasi-isomorphisms, it provides the abelian category of $1$-motives with torsion $\tM$. This definition has been generalized by Barbieri-Viale and Kahn in \cite{BVK} to any perfect field (inverting $p$-multiplication) and in \cite{BER} without inverting $p$-multiplication, but accepting connected finite $k$-groups in degree $-1$. Still in \cite{BVK} the authors introduced the category of $1$-motivic sheaves $\Shv_1$ and showed that $D^b(\Shv_1)$ is equivalent to $D^b(\tM)$ and, moreover, equivalent to the thick subcategory of Voevodsky's triangulated category of motives ${\rm DM}^\eff_{gm}(k)$ generated by motives of smooth curves. 
In \cite{BER} we extended the first equivalence to perfect fields of positive characteristic that are transcendental over their prime field.
In this paper we show that this equivalence result can be extended to the ``Laumon context'', \ie allowing additive groups and formal connected $k$-groups.  The interest in Laumon $1$-motives arises from the expectation that they will play an important role in a contravariant theory of motivic cohomology for singular varieties.
See \cite{BVB}, \cite{LEK}, for results in this direction.

%%%
\section{Generalized $1$-motives with torsion}\label{sec.deriv2}

Any formal (commutative) $k$-group $L$ may be viewed as an fppf sheaf on the category of affine $k$-schemes and hence it provides an fppf sheaf on the category of $k$-schemes. Let $L$ be a formal $k$-group with maximal connected subgroup $L^0$, and put $L^\et\df L/L^0$. Furthermore, let $L^\et_\tor$ be the maximal torsion subgroup of $L^\et$ and $L^\et_\fr:=L^\et/L^\et_\tor$.
We will extend constructions and results of \cite{BVK} and \cite{BER} to Laumon $1$\nobd-motives.

In the following by connected algebraic $k$-group we mean a connected, smooth, commutative $k$-group scheme.
\begin{defn}\label{def.geffma}
A \emph{generalized effective  $1$\nobd-motive with torsion} is a complex  of fppf sheaves on $\Spec(k)$, $M=[L\to G]$, where $L$ is a formal $k$-group and $G$ is a connected algebraic $k$-group. An \emph{effective morphism} is a map of complexes.
\end{defn}
These $1$-motives admit unipotent groups, \ie copies of $\G_a$, in the component of degree $0$ and formal connected $k$-groups, \ie powers of $\hat\G_a$, in the component of degree $-1$. See \cite{RUS} for a generalization. Denote $\tMaeff$ the category of (generalized) effective $1$\nobd-motives with torsion. The category of  Laumon $1$\nobd-motives $\Ma$ is the full subcategory of $\tMaeff$ consisting of those $[L\to G]$ with $L^\et$ torsion free (\cite{LAU}). Furthermore, if one requires $L$ to be \'etale and torsion free, and $G$ to be semiabelian, one gets the full subcategory of Deligne $1$\nobd-motives.

We now  introduce the category of $1$-motives with torsion. 
Recall first that an effective morphism $(f,g)\colon M\to M'$ of (generalized) effective $1$-motives with torsion  is a \emph{quasi-isomorphism (q.i.)} if  $f$ is a surjective map of formal $k$-groups, $g$ is an isogeny, and $\ker(f)=\ker(g)$ is a finite group scheme.

\begin{defn} The category of \emph{(generalized) $1$\nobd-motives with torsion} $\tMa$ is the localization of $\tMaeff$ with respect to the multiplicative class of quasi-isomorphisms.
\end{defn}
We have (see \cite{BVB}, 1.4.5):

\begin{thm}\label{thm.aca} \begin{itemize}
\item[(i)] $\tMa$ is an abelian category.
\item[(ii)] Any short exact sequence of $1$-motives in $\tMa$ can be represented by a sequence of generalized effective $1$-motives that is exact as sequence of complexes.
\item[(iii)] The canonical functor  $d\colon \Ma\to \tMa$ makes $\Ma$ an exact subcategory of $\tMa$.
\end{itemize}
\end{thm}

The category $\tMa$ was introduced in \cite{BVB} and shown to be equivalent to the category of Formal Hodge Structures of level $\leq 1$ when $k=\C$; see also \cite{FHS} for the ``torsion free'' case. This result  generalizes  equivalences in \cite{BVRS} and \cite{DEL} for Deligne's $1$\nobd-motives.

The proof of the above theorem uses  the properties of strict morphisms  (cf. \cite{BVK}, C.5.2, \cite{BER}, 1.10 for details).

\begin{defn}
An effective morphism $\varphi=(f,g)\colon M\to M'$ is \emph{strict} if $g$ has smooth connected kernel, \ie the kernel of $g$ is still a connected algebraic $k$-group.
\end{defn}

Let $(f,g)\colon M\to M'$ be an  effective map. Its kernel $[\ker^0(f)\to \ker^0(g)]$ is defined as the effective $1$-motive with torsion that is the identity component of $\ker(g)$ in degree $0$, and the pull-back of  $\ker^0(g)$ along $ \ker(f)\to\ker(g)$ in degree $-1$. The cokernel of $\varphi$ is the effective $1$-motive with torsion $[\coker(f)\to \coker(g)]$. One sees quite easily that the localization functor  ${}^t\!\M^{{\rm a}, \eff}\to \tMa$ preserves kernels. For cokernels this is true when $\varphi$ is strict. In general, the fact that any effective map factors as a strict morphism followed by a quasi-isomorphism (cf. \cite{BVRS}, 1.4) is often useful.

In order to give a nice description of the bounded derived category of  $\tMa$ we introduce the following category.

\begin{defn} Let $\M^{{\rm a}\star}$ denote the full subcategory of $\Ma$ whose objects are those $[u\colon L\to G]$ with $\ker (u)=0$.
\end{defn}
Observe that there are no non-trivial quasi-isomorphisms  in $\M^\star$.

\begin{lemma}\label{lem.maw}
$\M^{{\rm a}\star}$ and$ \Ma$ are generating subcategories of $\ \tMa$ closed under kernels and  closed under extensions.
\end{lemma}
\proof The categories  $\Ma$ and $\M^{{\rm a}\star}$ are clearly closed under kernels because, given a morphism $(f,g)\colon M_1\to M_2$ of Laumon $1$-motives,  the formal $k$-group $\ker^0(f)$ is the kernel of the composition of maps $\ker(f)\to\ker(g)\to \Phi$ where $\Phi$ the component group of $\ker(g)$. 
To check that the above categories are closed  under extensions one applies Theorem~\ref{thm.aca}.  The verification that they are generating subcategories requires more work.  Let $M=[u\colon {\hat\G_a}^r\times L^\et\to G]$ be a generalized $1$-motive with torsion. We need to show that there exists a $1$-motive $M'$ in $\M^{{\rm a}\star}$ and an epimorphism $M'\to M$ in $\tMa$. 
By replacing, if necessary,  $L^\et$ by a torsion free \'etale formal $k$-group that dominates it, we may assume that $L^\et$ is torsion free.   
Let  $v\colon L^\et\to T$ be an embedding with $T$ a torus (clearly such an embedding exists over a finite extension $k'$ of $k$, say $L^\et_{k'}=\Z^r\to \G_{m,k'}^r$ mapping $1\mapsto 1$ in each component, and then consider the induced homomorphism $L^\et\to T:= \Re_{k'/k}(  \G_{m,k'}^r)$  by restriction of scalars). Set $\iota\colon {\hat\G_a}^r\to \G_a^r$ the canonical embedding. Then the $1$-motive $M'=[u'\colon {\hat\G_a}^r\times L^\et\to \G_a^r\times T\times G]$, $u'(x,y)=(\iota(x),v(y),u(x,y))$,  is in $\M^{{\rm a}\star}$. Furthermore the map $(\id,p_G)\colon M'\to M$, with $p_G$ the projection on the third component, is a strict epimorphism of effective $1$-motives and hence it remains an epimorphism in $\tMa$. Therefore $\M^{{\rm a}\star}$  is a generating subcategory, and hence so too is $\Ma$. \qed

We will see in Section~\ref{sec.equiv} how $D^b(\tMa)$ can be described in terms of complexes in $\M^{{\rm a}\star}$.

 %%%%
\section{Generalized $1$-motivic sheaves}\label{sec.genmot}
 In this section we define the category of generalized $1$-motivic sheaves. It contains the category of $1$-motivic sheaves introduced in \cite{BVK}, \cite{BER}.  We will show in Section~\ref{sec.equiv} that it is an abelian category whose bounded derived category is equivalent to $D^b(\tMa)$.

\begin{defn}\label{def.motsha}
An fppf sheaf $\mathcal  F$ on $\Spec(k)$ is a \emph{generalized $1$-motivic sheaf} if (i)  there  exists a morphism of sheaves $b\colon G\to \mathcal F$ with $G$ a smooth connected  algebraic $k$-group and $E=\coker(b)$, $L=\ker(b)$ formal $k$-groups, (ii) the $2$-fold  exact sequence of sheaves
\begin{eqnarray}\label{eq.Fa}
\eta\colon \ 0\to L\by{a} G\by{b}\mathcal F\by{c} E\to 0 ,
\end{eqnarray}
becomes trivial when restricted (by pull-back) to $E^0$.

The morphism  $b$ is said to be \emph{normalized} if $\ker (b)$ is torsion free.
\end{defn}

Obviously condition (ii) above is superfluous if $\mathrm{Ext}^2(\hat\G_a,L)=0$  for any formal $k$-group $L$. Unfortunately we have no proof  for the vanishing of the above $2$-fold extension group of fppf sheaves.

Definition~\ref{def.motsha} is the natural analogue of the one in \cite{BVK}, \cite{BER}, except for the additional condition (ii). Recall that  an fppf sheaf over $k$ is $1$-motivic in the sense of \cite{BER} if and only if it is the cokernel of an injective map $F_1\to F_0$ with $F_1$ a formal \'etale $k$-group and $F_0$ an extension of a formal \'etale $k$-group by a semiabelian variety (cf. \cite{BER}, 4.7).  We will see in Proposition~\ref{pro.pres} that an fppf sheaf is generalized $1$-motivic if and only if it is the cokernel of an injective map $F_1\to F_0$ where $F_1$ is a formal $k$-group and $F_0$ is an extension of a formal $k$-group by a connected algebraic $k$-group. Hence the above definition is the correct analogue of those in  \cite{BVK}, \cite{BER}.  

Let $\Shv_1^{\rm a}$ denote the full subcategory of fppf sheaves on $\Spec(k)$ consisting of  (generalized) $1$\nobd-motivic sheaves and  $\Shv_0^{\rm a}$ the subcategory consisting of those $\mathcal F$ with $G=0$. The latter category is equivalent to ${\bf For}/k$.

\begin{examples}\label{exa}
$\bullet$ The category $\Shv_1^\fppf$ introduced in \cite{BER} is equivalent to the full subcategory of $\Shv_1^{\rm a}$ consisting of those $\mathcal F$ with $G$ semiabelian, $L$ and $E$ \'etale.

$\bullet$ Let $M=[L\to G]$ be a Deligne $1$\nobd-motive over $k$. It is shown in \cite{BER} that for  $M^\natural=[u\colon L\to G^\natural]$ a universal $\G_a$-extension of $M$ the sheaf $\coker(u)$ is the sheaf of $\natural$-extensions of $M$ by  $\G_m$ (\ie those extensions endowed with a connection compatible with the extension structure). Since $G^\natural$ is a connected algebraic $k$-group (but, in general, not a semiabelian 
variety), and $L$ is a discrete group, $\coker(u)$ is an example of a generalized
$1$\nobd-motivic sheaf that is not $1$\nobd-motivic in the sense of \cite{BER}, 2.1.

$\bullet$ Let $A$, $A'$ be dual abelian varieties over $k$ and $\hat A'$ the completion of $A'$ along the identity. The Laumon $1$\nobd-motive $[\hat A'\by{\iota} A']$, for $i$ the canonical embedding, is the Cartier dual of $[0\to A^\natural]$ where $A^\natural$ is the universal $\G_a$-extension of $A$ (\cf \cite{LAU}). The sheaf $\coker(\iota )$ is a generalized $1$\nobd-motivic sheaf. More generally, given a smooth connected algebraic $k$-group $G$, let $\hat G$ be its formal completion at the identity.  Then $[\iota \colon \hat G\to G]$ is a generalized $1$\nobd-motive and $\coker(\iota)$ is a $1$\nobd-motivic sheaf.
\end{examples}

Propositions 3.2.3 in \cite{BVK} and 2.3 in \cite{BER} generalize easily.

\begin{propose}\label{pro.motsha}
\begin{itemize}
\item[(i)] In Definition~\ref{def.motsha} we may choose $b$  normalized.
\item[(ii)] Given two $1$\nobd-motivic sheaves $\mathcal F_1, \mathcal F_2$, normalized morphisms $b_1\colon G_1\to \mathcal F_1$, $b_2\colon
G_2\to \mathcal F_2$ and a map $\varphi\colon \mathcal F_1\to
\mathcal F_2$ of fppf sheaves, there exists a unique homomorphism of
group schemes $\varphi_G\colon G_1\to G_2$ over $\varphi$. In particular a morphism in $\Shv_1^{\rm a}$ is uniquely determined by a morphism of the complexes \eqref{eq.Fa} with normalized $b$'s.
\item[(iii)] Given a $1$\nobd-motivic sheaf $\mathcal F$, a
 morphism  $b\colon G\to \mathcal F$  as above with $b$ normalized
is uniquely (up to isomorphisms) determined by $\mathcal F$.
\item[(iv)]  Let  $\varphi\colon \mathcal F_1\to \mathcal F_2$ be a morphism of fppf sheaves, where $\mathcal F_i$ fits into an exact sequence $\eta_i\colon 0\to L_i\to G_i\to \mathcal F_i\to E_i\to 0$, $i=1,2$, of the type \eqref{eq.Fa}.  Suppose that $\varphi$ induces an injective map $\varphi_L\colon L_1\to L_2$; if $\mathcal F_2$ is $1$-motivic, then so too is $\mathcal F_1$. Suppose that $\varphi_E\colon E_1\to E_2$, induced by $\varphi$, is surjective; if $\mathcal F_1$ is $1$-motivic, then the same holds for  $\mathcal F_2$.
\item[(v)]$\Shv_1^{\rm a}$ and $\Shv_0^{\rm a}$ are exact abelian subcategories of the category of fppf sheaves on $\Spec(k)$.
\end{itemize}\end{propose}
\proof
(i). If $b$ is not normalized, divide $G$ by $L^\et_\tor$. 
 (ii).  Consider the exact sequence \eqref{eq.Fa} for both $\mathcal F_i$, $i=1,2$. The composition $c_2\varphi b_1$ is trivial, hence the morphism $\varphi$ induces a morphism $\varphi_E\colon E_1\to E_2$ and a morphism $G_1\to G_2/L_2$. The latter lifts to a morphism $\varphi_G\colon G_1\to G_2$ because $\Ext^1(G_1,L_2)=0$ (\cf~\cite{BVB}, Lemma A.4.5). The uniqueness of  $\varphi_G$ follows from the vanishing of  $\Hom(G_1,L_2)$. 
(iii) follows from (ii).  
For (iv): suppose $\varphi_L$ injective and $\mathcal F_2$ $1$-motivic. Only property (ii) in Definition \ref{def.motsha} has to be checked for $\mathcal F_1$.  
Denote $\eta_i^0$ the pull-back of $\eta_i$ along $E_i^0\to E_i$. Then $\eta_1^0$ maps term by term to $\eta_2^0$.  Since $\mathcal F_2$ is $1$-motivic $\eta_2^0$ is trivial. Hence the same holds for the push-out of $\eta_1^0$ along $\varphi_L$. 
By the injectivity of $\varphi_L$ we deduce that $\eta_1^0$ is obtained from an extension of $E_1^0$ by the formal $k$-group $L_2/L_1$. Hence it is trivial. 
The proof of the other case is analogous.
 (v) is evident for $\Shv_0^{\rm a}$. Consider a morphism of generalized $1$-motivic sheaves $\varphi \colon \mathcal F_1\to
\mathcal F_2$. Define $G_3$ as the identity component of the kernel of $\varphi_G\colon G_1\to G_2$. We have a canonical map $b_3\colon G_3\to \ker(\varphi)=:\mathcal F_3$. Its kernel is still in $\mathcal{ CE}$. Define $G_4$ as the cokernel of $\varphi_G$. We have a canonical map $b_4\colon G_4\to \coker(\varphi)=:\mathcal F_4$. By a diagram chase one sees that both $\mathcal F_3$ and $\mathcal F_4$ fit into an exact sequence of the type \eqref{eq.Fa} and thanks to (iv) these sequences satisfy condition (ii) in Definition~\ref{def.motsha}.
\qed

We will show in the following paragraphs that any generalized $1$-motivic sheaf admits a
presentation $[ F_1\to F_0]$ where $F$ is a sheaf extension of a formal $k$-group by a smooth connected algebraic $k$-group and $F_1$ is a formal $k$-group. 
To prove this fact we will need several results on extensions.

\begin{lemma}\label{lem.extdd}
Let $\eta$ be an exact sequence of fppf sheaves over $\Spec(k)$ as in \eqref{eq.Fa}  with $L,E$ discrete. Then, there exists an epimorphism $\varphi\colon \tilde E\to E$ of discrete groups such that the pull-back  of $\eta$ along $\varphi$ is isomorphic to the  trivial extension.  
\end{lemma}
\proof \cite{BER}, Lemma 4.6.
\qed

In order to prove a similar result in the general case we need the following lemma.

\begin{lemma}\label{lem.H1H2}
 ${\mathrm H}^1(k,\hat \G_a)=0={\mathrm H}^2(k,\hat \G_a)$.
\end{lemma}
\proof Recall the spectral sequence
\begin{equation}\label{eq.spectral}
{\mathrm H}^p(X,\underline {\rm Ext}^q(Y,Z)) \Rightarrow {\rm Ext}^{p+q}_X(Y,Z). \end{equation}
Set $X=\Spec(k)$, $Y=\G_a$, $Z=\G_m$. Since $\underline{\rm Ext}^1(\G_a,\G_m)=0$ and ${\rm Ext}^{n}_k(\G_a,\G_m)=0$ for any $n$ (cf. \cite{BRE}), using the five terms exact sequence associated to \eqref{eq.spectral} (cf. \cite{EtC}, Appendix B)   one  proves the assertions. \qed

\begin{lemma} \label{lem.triv}
Let $\eta\colon 0\to L\to G\to {\mathcal F}\to E\to 0$ be an exact sequence of fppf sheaves over $\Spec(k)$ and $\mathcal F$ a generalized $1$-motivic sheaf. Then,  there exists an epimorphism $\varphi\colon \tilde E\to E$ of formal $k$-groups, $\tilde E$ torsion free, such that the pull-back of $\eta$ along $\varphi$ is isomorphic to the trivial extension. 
\end{lemma}
\proof
We may assume $E=E^\et$. Indeed, since $E=E^0\oplus E^\et$, there exist $1$-motivic sheaves \[\eta^0\colon 0\to L\to G\to {\mathcal F}^\et \to E^\et \to 0, \quad   \eta^\et\colon  0\to L\to G\to {\mathcal F}^0\to E^0\to 0\] such that $\mathcal F$, viewed as extension of $E$ by $\mathcal F^\star=G/L$, is the ``sum'' of $\mathcal F^0$ and $\mathcal F^\et$, \ie the push-out along the multiplication map $\mu\colon \mathcal F^\star\oplus \mathcal F^\star \to \mathcal F^\star$ of $\mathcal F^0\oplus \mathcal F^\et$. In particular $\eta$ is the ``sum'' of $\eta^0$ and $\eta^\et$, \ie the push-out along $\mu_L\colon L\oplus L\to L$ of $\eta^0\oplus \eta^\et$. By definition of $1$-motivic sheaves, $\eta^0$ is trivial.
It is then sufficient to prove the lemma for $E=E^\et$.

Let $f\colon \Spec(k')\to \Spec(k)$ be a Galois extension of degree $n$ such that $f^*E$ becomes constant.  The functor $f^*$ is exact (\cite{EtC}, II 2.6 (b)) and is the restriction functor; hence  $f^*\hat \G_a=\hat\G_{a,k'}$.
Enlarging $k'$, if necessary, we may assume that $\mathcal F(k')\to E(k')$ is surjective. Indeed ${\rm H}^1(k', f^*\mathcal F^\star)$ fits into a long exact sequence of cohomology groups, between   ${\rm H}^1(k', G_{k'})$ and ${\rm H}^2(k',f^*L)= {\rm H}^2(k',L^\et_{k'})$ (see Lemma~\ref{lem.H1H2}). Hence any element in ${\rm H}^1(k', f^*\mathcal F^\star)$ becomes trivial over a suitable finite Galois extension of $k'$; in particular, enlarging $k'$ we may assume that any element of a fixed finite set of generators of $E_{k'}$ becomes trivial in ${\rm H}^1(k', f^*\mathcal F^\star)$.  The surjectivity of $\mathcal F(k')\to E_{k'}$ implies the triviality of the exact sequence $0\to f^*\mathcal F^\star\to f^*\mathcal F\to E_{k'}\to 0$. Let $s\colon E_{k'}\to f^*\mathcal F$ be a section.  We may lift $s$ to a morphism  $s'\colon E_{k'}'\to f^*\mathcal F$ with $E_{k'}'\to E_{k'}$ a surjective morphism of formal $k$-groups with $E_{k'}' $ torsion free; simply choose $E_{k'}'=E_{k',\fr}\times \Z^r$ with $E_{k',\tor}=\oplus_{i=1}^r\Z/m_i\Z$. 
Consider now the trace map \[\tau\colon f_*f^*E\to E, \quad b\mapsto \sum_{\sigma\in Gal(k'/k)} \sigma(b) .\]
Observe that it is surjective because its base change to $k'$ is the map
\[ f^*f_*E_{k'}=\prod_{\sigma\in Gal(k'/k)} \!\!E_{k'} \to E_{k'}, \quad (a_\sigma)\mapsto \sum_\sigma a_\sigma.\]
The pull-back of $\eta$ along $f_*E_{k'}'\to f_*E_{k'}$ followed by $\tau$ is trivial because we have a map $f_*E_{k'}'\to f_*f^*\mathcal F\to \mathcal F$ where the first map is $f_*s'$ and the second map is the trace map. 
\qed

We prove now the main result on presentations of generalized $1$-motivic sheaves, thus generalizing   Proposition~4.7 of \cite{BER}.

\begin{propose}\label{pro.pres}
An fppf sheaf  $\mathcal F$ on  $\Spec(k)$ is $1$-motivic if and only if
\begin{eqnarray}\label{eq.motsha}
\mathcal F=\coker(F_1\by{u} F_0)
\end{eqnarray}
where $F_1$ is a formal $k$-group,  $F_0$ is an extension of a formal $k$-group by a smooth connected algebraic $k$-group $G$ and  $u$ is a monomorphism. We then have a diagram  with exact rows and columns
\begin{eqnarray}\label{dia.def}
\xymatrix{
0\ar[r]& L\ar[r]\ar@{^{(}->}[d]& G\ar[r]\ar@{.>}[dr]^b\ar@{^{(}->}[d]& \mathcal F^*\ar[r]\ar@{^{(}->}[d]& 0\\
 0\ar[r]  & F_1 \ar@{->>}[d] \ar[r]^u & F_0 \ar[r]\ar@{->>}[d]  & {\mathcal F}\ar[r]\ar@{->>}[d]&  0\\
0\ar[r]& F_1/L \ar[r] & F_0/G  \ar[r]& E\ar[r] & 0
 }\end{eqnarray}(\ref{dia.def}) with
$L,E$ formal $k$-groups. Moreover, we can always find a presentation
$F_1\to F_0$ of $\mathcal F$ with $F_1$ torsion free.
\end{propose}
\proof  To prove sufficiency, let $\mathcal F=F_0/F_1$ be as in the above diagram. 
Then $\mathcal F$ fits into an exact sequence
\begin{equation}\label{eq.Fa2}
\eta\colon \ 0\to L\to G\by{b}\mathcal F\to E\to 0 ,
\end{equation}
where $b$ is the composition $ G\to F_0\to \mathcal F$.
Furthermore, the existence of $F_0$ says that the pull-back of  the extension \eqref{eq.Fa2} along $F_0/G\to E$ is trivial. Applying Proposition~\ref{pro.motsha} (iv) one gets that the pull-back of \eqref{eq.Fa2} to $E^0$ vanishes.
Hence $\mathcal  F$ is $1$-motivic.

We now prove necessity of the condition. Let $\mathcal F$ be a $1$-motivic sheaf and consider the associated exact sequence \eqref{eq.Fa2}.  For the existence of a presentation as in \eqref{dia.def}, it is sufficient to see that, after pull-back along a suitable epimorphism $\psi\colon \tilde E\to E$ of formal $k$-groups, the $2$-fold extension in \eqref{eq.Fa2} is isomorphic to the  trivial one.  
Indeed let \[\psi^*(\eta)\colon 0\to L\to G\to\tilde{\mathcal F}\to \tilde E\to 0.\]
If the class of $\psi^*(\eta)$ is trivial, then $\tilde{\mathcal F}=\tilde F_0/L$ with $\tilde F_0$ extension of $\tilde E$ by $G$ and $\mathcal F$ is the quotient of $\tilde F_0$ by the formal $k$-group $\tilde F_0\times_{\tilde {\mathcal F}} \ker(\psi)$. 
To conclude the proof it is sufficient to recall that such a $\psi$ exists by Lemma~\ref{lem.triv}.

Observe futhermore that $\mathcal F=\coker(F_1/F_{1,\tor}\to F_0/F_{1,\tor})$.
\qed

We now introduce a useful subcategory of $\Shv_1^{\rm a}$. 

\begin{defn}
Let  $\Shv_1^{{\rm a}\star}$ be the full subcategory of $\Shv_1^{\rm a}$ consisting of those $1$-motivic sheaves such that $\coker(b)=0$, $b$ as in \eqref{eq.Fa}.
\end{defn}

Before proving that  the category $\Shv_1^{{\mathrm a}\star}$ has the ``dual'' properties of $\M^{{\mathrm a}\star}$,  we still need a preparatory result:

\begin{lemma}\label{lem.triv2}
Let $G$ be a connected algebraic $k$-group. Then $\mathrm{Ext}^1(\hat\G_a,G)=0$.
\end{lemma}
\proof
The cases $G=\G_a,\G_m$, were treated in \cite{BVB}.
Denote $\hat G$ the formal completion of $G$ along the unit section $\epsilon$.
We will see that any element in $\mathrm{Ext}^1(\hat\G_a,G)$ is obtained via push-out $\hat G\to G$ from an extension in $\mathrm{Ext}^1(\hat\G_a,\hat G)$.   Thus, since the latter group is trivial, this will imply the desired result.
Consider an exact sequence $0\to G\to H\to \hat \G_a\to 0$. Denote by  $\G_{a,n}$ the scheme $\Spec(k[x]/(x^{n+1}))$ and by $H_n$ the pull-back of $H$ along the canonical map $\G_{a,n}\to \hat \G_a$. This $H_n$ is a scheme because $\G_{a,n}$ has dimension $0$ (\cite{EtC}, III, 4.3).  Furthermore $H_n$ is smooth over $\G_{a,n}$ (\cite{RAY}, VI, 1.2). Clearly $G=H_0$  over $\Spec(k)$ admits a section $s_0$, the unit section. By the smoothness of any $H_n$ over $\G_{a,n}$,  this section lifts to a compatible tower of sections $s_n\colon  \G_{a,n}\to H_n$, $n\geq 0$.  Define $\gamma_n\colon \G_{a,n}^2\to G$ as
\[\gamma_n(a,b)= s_{2n}(a+b)-s_n(a)-s_n(b),\]
where the sum $a+b$ lives in $\G_{a, 2n}$. It satisfies $s_n(0)=\epsilon\in G(k)$. 
In this way we define a map $\gamma\colon \hat \G_a\to G$ of contravariant functors of sets satisfying the usual properties of a factor set.  The map $\gamma$ factors through $\hat G$. Indeed, if $U=\Spec(A)$ is an open affine neighborhood of the unit section $\epsilon$ of $G$, then $\gamma$ factors through $U$. Let $m$ be the maximal ideal of $A$ corresponding to $\epsilon$.  One has $\gamma_n^*(m)\subseteq (x)$ and $\gamma_n$ factors then through $\Spec(A/m^n)$. 
 Hence $\gamma$ provides an extension $\eta$ of $\hat G$ by $\hat \G_a$, \ie  the extension we started with is the push-out along $\hat G\to G$  of $\eta$.  Since $\eta$ is trivial, the lemma follows.
\qed

\begin{lemma}\label{lem.shva10}
$\Shv_1^{{\mathrm a}\star}$  is a cogenerating subcategory of $\Shv_1^{\rm a}$, closed under cokernels and  closed under extensions.
\end{lemma}
\proof The only non trivial fact is  that  $\Shv_1^{{\rm a}\star}$ is cogenerating, \ie for any $1$\nobd-motivic sheaf $\mathcal F$ there exists a $\mathcal F'$ in $\Shv_1^{{\rm a}\star}$ and a monomorphism $\varphi\colon \mathcal F\to \mathcal F'$. To prove this, consider the diagram \eqref{dia.def} associated to $\mathcal F$.  It is sufficient to prove that $F_0$ embeds in a $1$-motivic sheaf $\mathcal F_0 $ in $\Shv_1^{{\rm a}\star}$ and then set $\mathcal F'=\mathcal F_0/F_1$. So assume $\mathcal F=F_0, L=0$.
 If  $E=E^0\cong \hat \G_a^s$, then $\mathcal F=F_0$ is the product $G\times E^0$ by Lemma~\ref{lem.triv2}. Now embed $F^0$ into $G\times \G_a^s$. If $E=E^\et$ we proceed as in the proof of \cite{BER}, 2.8.
Suppose $E_\tor$ is $n$-torsion. Denote by ${}_nF_0$ the kernel of the $n$-multiplication on $F_0$. It is an extension of $E_\tor $ by ${}_n G$ since the $n$-multiplication is surjective on $G$. Furthermore $F_0/{}_nF_0=G$. For any embedding $j\colon {}_nF_0\to B$ into an abelian variety, the push-out of $F_0$ along $j$ is an extension $\tilde F$ of $G$ by $B$, hence a connected algebraic $k$-group.
Suppose now that $E$ is torsion free \'etale. Over a suitable finite galois extension $k'$ of $k$ the group $E$ becomes isomorphic to $\Z^r$ and $F_{0,k'}\cong G_{k'}\times E_{k'}$. Embed  $E_{k'}$ into an abelian variety $B'_{k'}$. Then $F_{0 ,k'}$ embeds into the connected algebraic $k'$-group $H_{k'}=G_{k'}\times B_{k'}$. As a consequence $F_0$ embeds into the Weil restriction $\Re_{k'/k}(H_{k'})$ that is  a connected algebraic $k$-group. 
The general case is obtained by devissage from the above cases: first one considers the pull-back $F_0'$ of $F_0$ along $E^0\to E$. We have seen that  $F_0'$ embeds into a connected algebraic $k$-group $H_1$ and hence $F_0$ embeds into an extension $F_0''$  of $E^\et$ by $H_1$. Consider now the pull-back of $F_0''$ along $E^\et_\tor\to E^\et$. We have seen that $F_0''$ embeds into  a connected algebraic $k$-group $H_2$
and hence $F_0''$ embeds into an extension $F_0'''$ of $E^\et_\fr$ by $H_2$. Finally $F_0'''$ embeds into a connected algebraic $k$-group $H_3$ and the embedding we are looking for is the composition $F_0\to F_0''\to F_0'''\to H_3$. 
\qed

It follows from \eqref{eq.Fa} that any generalized $1$\nobd-motivic sheaf $\mathcal F$ can be viewed as an extension of a formal $k$-group $E$ by a $1$\nobd-motivic sheaf $\mathcal F^\star$ in $\Shv_1^{{\mathrm a}\star}$.
The category $\Shv_1^{{\mathrm a}\star}$ is interesting for us because of the following fact: 

\begin{lemma}
The functor ${\rm H}_0$ mapping a generalized effective $1$-motive
with torsion $[L\to G]$ to the generalized $1$-motivic sheaf $G/L$
provides an equivalence between $\M^{{\rm a}\star}$ and the full
subcategory $\Shv_1^{{\rm a}\star}$ of $\Shv_1^{\rm a}$.
\end{lemma}
\proof It follows from Proposition~\ref{pro.motsha} that
${\mathrm H}_0\colon \M^{{\rm a}\star} \to \Shv_1^{{\rm a}\star}$ is fully faithful.
Now, given a $1$-motivic sheaf
$\mathcal F$ in $\Shv_1^{{\rm a}\star}$ and a normalized morphism $b\colon
G\to \mathcal F$ as in \eqref{eq.Fa}, the $1$-motive
$[u\colon \ker (b)\to G]$ satisfies $\ker (u)=0$ and $\coker
(u)=\mathcal F$. Hence ${\rm H}_0$ is essentially surjective. \qed

%%%
\section{Equivalence on bounded derived categories}\label{sec.equiv}
We can now generalize the results in \cite{BER} on the equivalence of bounded derived categories to Laumon $1$-motives. 

\begin{lemma}\label{lem.der1a}
Let $N^b(\Shv_1^{{\rm a}\star})$ be the full subcategory of $K^b(\Shv_1^{{\rm a}\star})$ consisting
of complexes that are acyclic as complexes of $1$-motivic sheaves.
 The natural functor \[K^b(\Shv_1^{{\rm a}\star})/N^b(\Shv_1^{{\rm a}\star})\to D^b(\Shv_1^{\rm a})\]
 is an equivalence of categories.
\end{lemma}
\proof
The assertion follows immediately from Lemma~\ref{lem.shva10} and \cite{KS}, Lemma 13.2.2.
\qed

 \begin{lemma}\label{lem.der2a}
 Let $N^b( \M^{{\rm a}\star})$ denote the full subcategory of
$K^b( \M^{{\rm a}\star})$ consisting of complexes that are acyclic as complexes of
generalised $1$-motives with torsion.
The natural functor \[K^b( \M^{{\rm a}\star})/N^b( \M^{{\rm a}\star})\to D^b(\tMa)\] is an equivalence of
 categories.
\end{lemma}
\proof
This proof is ``dual'' to the previous one.  The``dual" conditions required in \cite{KS}, 13.2.2 ii), are satisfied thanks to Lemma~\ref{lem.maw}. 
Furthermore, the ``dual" statement of \cite{KS}, 13.2.1, holds
as well, and finally one applies \cite{KS}, 10.2.7 ii).
\qed

Since  Lemma~\ref{lem.maw} applies to $\Ma$ as well, we can replace $ \M^{{\rm a}\star}$ by $\Ma$ in the previous proof. 

\begin{lemma}\label{lem.der3a}
Let $N^b(\Ma)$ be the full subcategory of $K^b(\Ma)$ consisting of complexes that are acyclic as complexes of generalized $1$-motives with
torsion. The natural functor \[K^b(\Ma)/N^b(\Ma)\to D^b(\tMa)\] is an
equivalence of categories.
\end{lemma}

To conclude the comparison between the bounded derived categories of $\tMa$ and $\Shv_1^{\rm a}$
there remains the verification that the functor ${\rm H}_0$, which maps a generalized effective $1$-motive with torsion $[L\to G]$ to the generalized $1$-motivic sheaf $G/L$, 
provides an equivalence between $\M^{{\rm a}\star}$ and the full
subcategory $\Shv_1^{{\rm a}\star}$ of $\Shv_1^{\rm a}$.

\begin{lemma}\label{lem.acya}
 Let $M^\bullet$ be a complex in $K^b(\M^{{\rm a}\star})$.
Then $M^\bullet \in N^b(\M^{{\rm a}\star})$ if and only if ${\mathrm
H}_0(M^\bullet)\in N^b(\Shv_1^{{\rm a}\star})$. In particular ${\mathrm
H}_0$ induces an equivalence of categories
\[ K^b(\M^{{\rm a}\star})/N^b(\M^{{\rm a}\star})\to K^b(\Shv_1^{{\rm a}\star})/N^b(\Shv_1^{{\rm a}\star}).\]
\end{lemma}
\proof
This proof follows \cite{BER}, 3.5, with minor changes.
Let $M^\bullet\colon (\cdots\to M^ i\stackrel{d^i}{\to} M^{i+1}\to\cdots)$ be a complex in $K^b(   \M^{{\rm a}\star})$. Observe that $\ker (d^i)\in    \M^{{\rm a}\star}$ for any $i$ by Lemma~\ref{lem.maw}.
If $M^\bullet$ is acyclic in $K^b(\tMa)$, up to isomorphisms, $\coker (d^i)$ is in $   \M^{{\rm a}\star}$,
where $\coker (d^i)$ denotes here the cokernel of $d^i$ in $\tM$.
Hence, in order to prove that $H_0(M^\bullet)$ is an acyclic complex in $\Shv_1^{\rm a}$ it is sufficient to check the case of a short exact sequence $M^\bullet$ where $M^i$, $i=1,2,3$, are all in $\M^{{\rm a}\star}$. Furthermore the maps $d^i$ are all effective because $\M^{{\rm a}\star}$ is a full subcategory of $\tMa$.
By Theorem~\ref{thm.aca} (ii) the sequence $M^\bullet$ is represented up to q.i. by a sequence $0\to {\tilde M}^0\by{d^0}{\tilde M}^1\by{d^1} {\tilde M}^2\to 0$ of effective $1$-motives with torsion ${\tilde M}^i=[\tilde u_i\colon \tilde L_i\to \tilde G_i]$ that is exact as sequence of complexes; moreover $\ker(\tilde u_i)=0$ because ${\tilde M}^i$ is q.i. to $M^i$  and $\ker(u_i)=0$. 
The sequence  ${\mathrm H}_0(M^\bullet)$ is just the sequence $0\to \coker(\tilde u_0)\to \coker(\tilde u_1)\to \coker(\tilde u_2)\to 0$ and the latter is exact because of the usual ker-coker sequence.

For the converse,  let $M^\bullet\in K^{b}(   \M^{{\rm a}\star})$ and suppose that $H_0(M^\bullet)$ is acyclic in $K^b(\Shv_1^{{\rm a}\star})$. We prove that $M^\bullet$ is acyclic by induction on the length of the complex.  Suppose first that the length is less or equal $2$, \ie  
 \[M^\bullet\colon
0\to {M}^1\by{d^{1}}{M}^{2}\by{d^{2}} {M}^3\to 0.\]
Assuming the exactness of $H_0(M^\bullet)$, one checks directly that $d^1$ is a strict monomorphism (more precisely, the kernel of $G_1\to G_2$ is trivial as well as the kernel of $L_1\to L_2$). Furthermore the cokernel $\bar M^2$ of $d^1$ in $\tM^{{\rm a},\eff}$ (or equivalently in $\tMa$) satisfies $\ker(\bar u_2)=0$ and $H_0(\bar M^2)=H_0(M^3)$; hence it is q.i.  to \ to $M^3$ by Proposition~\ref{pro.motsha} (iii). Hence $M^\bullet$ is exact in $\tMa$.

Suppose now that the result is true for complexes of length $n\geq 2$ and assume $M^\bullet$ has length $n+1$.  Since the  category $\Shv_1^{{\rm a}\star}$ is closed under cokernels, the sheaf 
\[{\mathcal K}^{i}:=\ker({\mathrm H}_0(d^{i}))=\coker({\mathrm H}_0(d^{i-2}))\]
 is in $ \Shv_1^{{\rm a}\star}$. By  Proposition \ref{pro.motsha} there exists a unique (up to isomorphism) $1$-motive $\tilde M^i$  in $   \M^{{\rm a}\star}$ such that ${\mathrm H}_0(\tilde M^i)=K^i$ and the morphism $d^{i-1}$ factors through $\tilde M^i$ by Proposition \ref{pro.motsha} (ii).   Clearly ${\mathcal K}^{i}=0$ and $\tilde M^2=M^1$.
We then have the following complexes
 \[0\to \tilde M^{1}\to M^{2}\to \tilde M^{3}\to 0, 
\quad 0\to \tilde M^{3}\to M^{4}\to \dots \to M^{n+2}\to 0
\] 
both inducing exact sequences of $1$-motivic sheaves via $H_0$. By the induction hypothesis they are acyclic in $K^b(   \M^{{\rm a}\star})$ and hence so too is $M^\bullet$. 
\qed

We can now prove the main result of the paper, thus generalizing
\cite{BER}, Theorem~3.6, to Laumon $1$-motives:

\begin{thm}\label{thm.maina}
Let $N^b(\Ma)$ be the subcategory of $K^b(\Ma)$ consisting of complexes that are acyclic as complexes in $K^b(\tMa)$. Denote by $D^b(\Ma)$ the localization $K^b(\Ma)/N^b(\Ma)$.  One then has the following equivalences of categories
\[D^b(\Ma)\cong D^b(\tMa)\cong  D^b(\Shv_1^{\rm a}).\]
\end{thm}
\proof Thanks to Lemmas \ref{lem.der1a}, \ref{lem.der2a}, \ref{lem.der3a}, \ref{lem.maw} and~\ref{lem.acya}, we have
\begin{align*}D^b(\Ma):= K^b(\Ma)/N^b(\Ma)&\cong D^b(\tMa)\cong
  K^b(\M^{{\rm a}\star})/N^b(\M^{{\rm a}\star})\\ &\cong
K^b(\Shv_1^{\rm a\star})/N^b(\Shv_1^{\rm a\star})\cong
 D^b(\Shv_1^{\rm a}).
 \end{align*}
\qed

{\bf Acknowledgments:}
The author would like to thank L. Barbieri-Viale for directing her attention to this subject and for helpful discussions. She also thanks P. Jossen for pointing out a mistake in a previous version of the paper.

 \end{document}